\documentclass[11pt,amstex,amscd]{amsart}

\newcommand{\ra}{\rightarrow}

\newcommand{\CP}{\mathbb{CP}}

\newcommand{\BZ}{\mathbb{Z}}

\newcommand{\Spinc}{\ensuremath{\mathrm{Spin}^c}}

\newcommand{\Pic}{\mbox{Pic}}

\newcommand{\T}{{\ensuremath{\textsc{T}}}}

\usepackage{latexsym}
\usepackage{amsmath}
\usepackage{amscd}
\usepackage[all]{xy} 
\usepackage{amsfonts}
\usepackage{amssymb}
\usepackage{amsbsy}
\usepackage{graphicx}

\newtheorem{theorem}{Theorem}
\newtheorem{lemma}[theorem]{Lemma}
\newtheorem{corollary}[theorem]{Corollary}

\newtheorem{conjecture}[theorem]{Conjecture}

\newcounter{results}\stepcounter{results}

\begin{document}

\DeclareGraphicsExtensions{.eps}


\title{Seiberg--Witten vanishing theorem for $S^1$--manifolds with fixed points}

\author{Scott Baldridge}
\address{Department of Mathematics,  Indiana University,
 Bloomington, Indiana 47405, USA }
\email{\rm{sbaldrid@indiana.edu}}
\date{2001}

\maketitle

\begin{abstract}In this paper we show that the Seiberg--Witten
invariant is zero for all smooth  4--manifolds with $b_+{>}1$
which admit circle actions that have at least one fixed point.
Furthermore, we show that all symplectic 4--manifolds which admit
circle actions with fixed points are rational or ruled, and thus
admit a symplectic circle action.
\end{abstract}


\section{Introduction}

This paper addresses two problems in 4--manifold theory.  The
first concerns the computation of 4--dimensional diffeomorphism
invariants. Ever since the introduction of Donaldson invariants in
the early 1980's, efforts to calculate diffeomorphism invariants
centered upon large classes of smooth manifolds that have some
additional structure. One such class of manifolds thought to have
promise was 4--manifolds with effective circle actions, but the
extra structure given by such manifolds turned out to be
insufficient for calculating Donaldson invariants.

With the introduction of Seiberg--Witten invariants in 1994, old
problems were revisited with new hope.   {S.~Donaldson} showed how
to calculate the Seiberg--Witten invariants in the simplest case
where the 4--manifold was a product of a 3--manifold and a circle
\cite{sw:swand4man}. In 1997 T. Mrowka, P. Ozvath, and B. Yu
\cite{sw:sw_inv_seifert_space}, and simultaneously L. Nicolaescu
\cite{sw:abiabatic_limits}, studied the 3--dimensional
Seiberg--Witten equations of Seifert--fibered spaces.  In 2001,
Seiberg--Witten invariants for manifolds with free and fixed point
free circle actions  were calculated in \cite{sw:circleactions},
\cite{sw:sworbcirc}. In this paper we finish this line of research
by calculating the Seiberg--Witten invariants for $S^1$--manifolds
with fixed points. We prove:

\begin{theorem}
If X is a smooth, closed, oriented 4--manifold with $b_+{>}1$ and
admits a circle action which has at least one fixed point, then
the Seiberg--Witten invariant vanishes for all $\mbox{Spin}^c$
structures.\label{thm:1}
\end{theorem}

\medskip

When the action on $X^4$ is free, the quotient of the
$S^1$--action is a smooth, closed 3--manifold $Y$.  In this case
$X$ can be thought of as a unit circle bundle of a complex line
bundle over $Y$.  Hence $X$ with the given circle action is
specified by the quotient $Y$ and the Euler class $\chi \in
H^2(Y;\BZ)$ of the line bundle.

A fixed point free circle action will have nontrivial isotropy
groups, forcing the orbit space to be an orbifold rather than a
manifold. As in the free case,  a manifold with a fixed point free
$S^1$--action can still be considered a unit circle bundle, but it
is furthermore a principal $S^1$--bundle of an orbifold line
bundle over a $3$--dimensional orbifold.  In this setup,
$H^2(Y;\BZ)$ is replaced by the group $\Pic^t(Y)$ which records
local data around the singular set (see \cite{sw:sworbcirc}).  The
manifold $X$ with given $S^1$--action is then determined by the
orbifold $Y$ and the Euler class of the orbifold line bundle which
is now an element of $\Pic^t(Y)$.

Manifolds which admit circle actions with fixed points come with
more complicated local data, yet this extra structure gives more
control when computing the Seiberg--Witten invariants. The
intuition behind the calculation is to find an essential sphere of
nonnegative square; the existence of such spheres force the
Seiberg--Witten invariants to vanish.

Theorem \ref{thm:1} combined with the formula derived in
\cite{sw:sworbcirc} gives the means for calculating the
Seiberg--Witten invariants of {\em any} $S^1$--manifold $X$ with
$b_+{>}1$.  See \cite{sw:4man_kirby_calc} for an introduction to
Seiberg--Witten theory or \cite{sw:notes_on_sw_theory} for a more
detailed analysis.

\begin{theorem}[General Formula]
Let $X$ be a  smooth, closed, oriented $4$--manifold with
$b_+{>}1$ and a smooth circle action.

\begin{enumerate}
\item If the action has fixed points, then $SW_X(\xi)=0$ for any
\Spinc\ structure $\xi$.

\item If $X$ has a fixed point free $S^1$--action, let $Y^3$ be
the orbifold quotient space and suppose that $\chi \in \Pic^t(Y)$
is the orbifold Euler class of the circle action. If $\xi$ is a
\Spinc\ structure over $X$ with $SW^4_X(\xi)\not=0$, then
$\xi=\pi^*(\xi_0)$ for some \Spinc\ structure on $Y$ and
\begin{eqnarray*}
SW^4_X(\xi) = \sum_{\xi' \equiv \xi_0 \mod \chi} SW^3_Y(\xi'),
\end{eqnarray*}
where $\xi'-\xi_0$ is a well--defined element of $\Pic^t(Y)$. See
\cite{sw:sworbcirc} with respect to the $b_+{=}1$ fixed point free
case.
\end{enumerate}\label{thm:general}
\end{theorem}

\medskip

A theorem of Taubes \cite{sw:sw_inv_and_symp_form} says that the
\Spinc\ structure associated with the first Chern class of a
symplectic 4--manifold must have Seiberg--Witten invariant $\pm
1$. Symplectic 4--manifolds always have $b_+{>}0$ because the
wedge product of symplectic form with itself is the volume form.
Putting these facts together with Theorem \ref{thm:1} implies the
following corollary.

\begin{corollary} A symplectic 4--manifold which admits a circle
action with fixed points must have $b_+{=}1$. \label{cor:b_+=1}
\end{corollary}

This application shows the usefulness of the Seiberg--Witten
invariants  even when they vanish.  It also provides a nice
introduction to the second problem: classifying symplectic
4--manifolds which admit circle actions with fixed points.

Such manifolds have been extensively studied. If the circle action
preserves the symplectic form in the sense that the generating
vector field of the action $\T$ satisfies $\mathcal{L}_\T \omega =
d\iota_{\T} \omega =  0$ then the action is called {\em
symplectic}.  If in addition, $\iota_\T \omega$ is exact, then the
action is called {\em Hamiltonian} because there is an $f \in
C^\infty(X)$ such that $\iota_\T \omega = df$.  {D.~McDuff} showed
for 4--manifolds that the existence of a fixed point implied that
a symplectic circle action must also be a Hamiltonian action
\cite{symp:momentmap}.  The same is true in any dimension if the
action is semifree with only isolated fixed points
\cite{symp:circact:semifree_symp_circ_act}.

In 1990 M. Audin classified 4--manifolds with symplectic circle
actions having fixed points using McDuff's result and symplectic
reduction techniques \cite{symp:circact:4man_ham_symp} (or
\cite{symp:circact:4man_admit_mom_map}). But it remained unknown
whether every symplectic 4--manifold with a circle action having
fixed points also admitted a symplectic circle action. Corollary
\ref{cor:b_+=1} attacks this classification problem from the
opposite direction by restricting which $S^1$--manifolds can have
symplectic structures. It turns out that the restrictive theory
can be made to meet with Audin's theory and thus gives a complete
classification.

Theorem \ref{thm:1} can be profitably restated as saying that
every manifold with $b_+{>}0$ which admits a circle action with
fixed points has an essential sphere of nonnegative square. A
theorem of T.J. Li \cite{sw:symp:embedded_spheres} implies that
symplectic 4--manifolds with this condition must be rational or
ruled. Thus:

\begin{theorem}Every symplectic 4--manifold which admits a circle
action with at least one fixed point is $\CP^2$, an $S^2$--bundle
over a surface, or $\overline{\CP}^2$ blowups of $\CP^2$ or an
$S^2$--bundle over a surface.\label{thm:3}
\end{theorem}

While the circle action in Theorem \ref{thm:3} is not necessarily
symplectic, the manifolds listed above are exactly the ones which
Audin classified (c.f. \cite{symp:circact:4man_ham_symp} or
\cite{symp:circact:top_of_torus_act}). Putting the two ideas
together gives:

\begin{corollary}
Every symplectic 4--manifold which admits a circle action with at
least one fixed point also admits a symplectic circle
action.\label{cor:symp_circ_act}
\end{corollary}

This corollary leads us to an interesting question posed by C.
Taubes \cite{sw:geometry_of_sw}: if $Y^3\times S^1$ is symplectic,
does $Y$ fiber over the circle?  Partial positive results have
been posted in \cite{symp:symp_lef_fib_S1xM},
\cite{symp:lefschetz_cpx_strct_seifert}. One can ask a much more
general question for any $S^1$--manifold (c.f.
\cite{sw:circleactions}):

\begin{conjecture} Every symplectic 4--manifold which admits a
circle action also admits (possibly a different) symplectic form
and circle action which are symplectic with respect to each other.
\end{conjecture}

If the conjecture holds it effectively answers Taubes' question.
Corollary \ref{cor:symp_circ_act} proves the conjecture when the
action has fixed points.  While the full proof of this conjecture
seems out of reach at the moment, Corollary
\ref{cor:symp_circ_act} does lend support to the hypothesis that
symplectic 4--manifolds with
$S^1$--actions are very special.\\

\medskip

\noindent {\em Acknowledgements.} I would like to thank Ron
Fintushel for his support through the summer of 2001 and for his
many helpful discussions during that time.  I would also like to
thank Tian--Jun Li for reading an earlier draft of this paper and
pointing out his results on embedded spheres could possibly be
combined with the proof of Theorem \ref{thm:1} to get Theorem
\ref{thm:3}.

\section{Proofs}

The theorem below proves both Theorem \ref{thm:1} and Theorem
\ref{thm:3}:

\begin{theorem}Suppose $X$ is a smooth, closed, oriented, 4--manifold with $b_+(X){>}0$ which admits
a smooth $S^1$--action with at least one fixed point.  Then $X$
contains an essential embedded sphere of nonnegative
self--intersection.\label{lem:main}
\end{theorem}

If a 4--manifold $X$ has an essential embedded sphere of
nonnegative self--intersection and $b_+(X){>} 1$, then the
Seiberg--Witten invariant vanishes (c.f. \cite{sw:class_turkey}),
proving Theorem \ref{thm:1}. Likewise, the existence of such a
sphere in a symplectic 4--manifold implies that the manifold is
rational or ruled \cite{sw:symp:embedded_spheres}.

The proof of Theorem \ref{lem:main} relies heavily upon R.
Fintushel's foundational work on 4--manifolds with $S^1$--actions
\cite{circact:Circ_act_on_four_man},
\cite{circact:Class_of_circ_4man}. Every such 4--manifold $\pi:X
\ra Y$ has a quotient 3--manifold $Y$ together with the following
data (altogether called a legally weighted 3--manifold $Y$):

\begin{enumerate}
\item  A finite collection of weighted arcs and circles in
$\mbox{int } Y$.

\item A finite set of isolated fixed points in $\mbox{int } Y$ disjoint from the sets of
(1).

\item A class (the ``Euler class'') $\chi \in H_1(Y, S)$ where $S$ is
the union of $\partial Y$, points of (2), and arcs of (1).

\end{enumerate}

The endpoints of weighted arcs are fixed points and there may be a
finite set of fixed points in the interior of a weighted arc or
circle.   To each component of an arc or circle that runs between
two fixed points (without passing through another fixed point) one
assigns a weight which records the local isotropy data of the
circle action.  If a weighted circle does not contain any fixed
points, it is called {\em simply--weighted}; otherwise it is {\em
multiply--weighted}.

To each component of (1) and (2) there is an associated index
which is simply the Euler number of the principal $S^1$--bundle
over the boundary of a tubular neighborhood of the component.
Similarly, an index can be given to each boundary component (see
Section 3.3 of \cite{circact:Circ_act_on_four_man}).  The sum of
the indices of all of the components is equal to 0 \cite[Page
380]{circact:Class_of_circ_4man}.

We also need Fintushel's classification result for simply
connected 4--manifolds with circle actions
\cite{circact:Circ_act_on_four_man},
\cite{circact:Class_of_circ_4man}.

\begin{theorem}[Fintushel] Let $S^1$ act smoothly on a simply
connected 4--manifold $X$, and suppose the quotient space $Y
\simeq S^3 $ is not a counterexample to the 3--dimensional
Poincar\'{e} conjecture. Then $X$ is a connected sum of copies of
$S^4$, $\CP^2$, $\overline{\CP}^2$, and $S^2\times S^2$.
\label{thm:Fintushel}
\end{theorem}

We can already use this result to prove Theorem \ref{thm:1} when
$X$ is simply connected. In that case, $X$ is the equivariant
connect sum of two $b_+>0$ pieces where one of the quotient spaces
is $S^3$. Thus there is a $\CP^2$ or $S^2\times S^2$ summand in
that piece with an essential embedded sphere of nonnegative
square, as required. We are now ready to generalize this idea.

\begin{proof}
Let $X$ be a smooth, closed 4--manifold with a smooth
$S^1$--action whose quotient is a legally weighted 3--manifold
$Y$. First we show how to reduce to the case where there are no
multiply--weighted circles in $Y$.

{Suppose $Y$ contains a multiply--weighted circle $C$ with 3 or
more fixed points.  Note that $C$ could represent a nontrivial
class in $H_1(Y;\BZ)$ or be embedded in $Y$ nontrivially as in
Figure \ref{fig:weightedorbit1}.

\begin{figure}[h]
\begin{center}
   \includegraphics{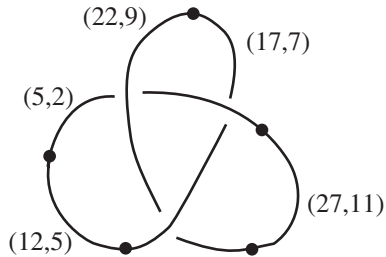}
   \caption{Example of a multiply--weighted circle}
   \label{fig:weightedorbit1}
\end{center}
\end{figure}
}

In this situation $X$ can be decomposed into an equivariant
connect sum  of two 4--manifolds $X=X_0 \# N_1$, both with circle
actions. The weighted orbit space of $X_0$ is the same as before
except the weighted circle $C$ has exactly two fixed points; the
weighted orbit space of $N_1$ is $S^3$ with a trivially embedded
multiply--weighted circle with the original weights (Figure
\ref{fig:weightedorbit2}).  The equivariant connect sum is
performed by removing the preimage of a $D^3$ neighborhood of the
fixed point between the weights $(5,2)$ and $(12,5)$ from both
$X_0$ and $N_1$ and then gluing equivariantly along the $S^3$
boundary.

\noindent\begin{figure}[h]
\begin{center}
   \includegraphics{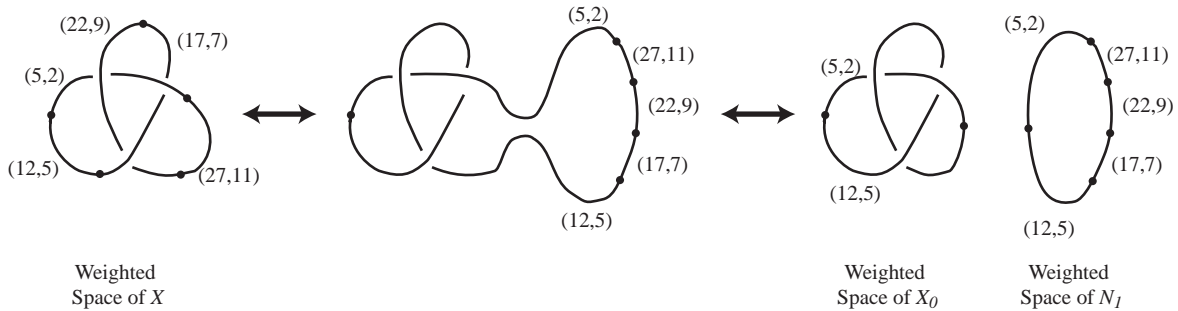}
   \caption{Decomposing $X$ into $X_0 \# N_1$}
   \label{fig:weightedorbit2}
\end{center}
\end{figure}

If we repeat this for all multiply--weighted circles in $X$ with 3
or more fixed points, we can decompose $X$ into

\[X= X_0 \# N_1 \# N_2 \cdots \# N_k\]

\noindent and we can further assume each $N_i$ decomposes into
connect sums of $\CP^2$, $\overline{\CP}^2$, $S^2 \times S^2$, or
$S^4$ by Theorem \ref{thm:Fintushel}. If any of the $N_i$'s have a
$\CP^2$ or $S^2\times S^2$ factor, then $X$ has an essential
embedded sphere of nonnegative square.  So we may assume that $X$
decomposes into $X_0 \# k\overline{CP}^2$ where $b_+(X_0)>0$ and
$k\geq 0$.  If there is an essential sphere of nonnegative square
in $X_0$, then such a sphere survives in any of its blowups. Hence
it is enough to find a sphere of nonnegative square in $X_0$.

We can get rid of the  multiply--weighted circle in $X_0$ using P.
Pao's ``replacement trick'' \cite{circact:nonlincircact}. He
noticed that the same 4--manifold often admits many different
$S^1$--actions and showed how to replace a complicated circle
action with a simpler one.

\begin{lemma}[Pao] Let $X$ be a 4--manifold with an $S^1$--action
whose weighted orbit space $Y$ contains a weighted circle $C$ with
exactly two fixed points.  Then $X$ admits a different
$S^1$--action whose weighted orbit spaces is either $Y$ with $C$
replaced with a pair of isolated fixed points or $Y\setminus
\mbox{int } D^3$ with C removed from the collection of weighted
circles. \label{lem:Pao}
\end{lemma}

Thus we can work with $X_0$, rename it $X$, and assume it has only
boundary components, weighted arcs, isolated fixed points, or
simply--weighted circles.  If there are two or more boundary
components in the quotient, the preimage of an arc that runs from
one boundary component to another boundary component is an
essential sphere of square zero; so we may reduce to the case
where $Y$ has only one boundary component. This case can be
eliminated using a short argument in the lemma below.

Thus we only need to consider $S^1$--manifolds $X$ with weighted
arcs, isolated fixed points, and simply--weighted circles. Isotope
the arcs and fixed points into a smooth ball $D^3 \subset Y$ and
enclose them by a sphere $\partial D^3 \subset Y$. Since the sum
of the indices of the components is 0, the Euler class of the
preimage of the sphere in $X$ is zero, i.e., we can realize $X$ as
the fiber sum of two manifolds $X_1$ and $N$ by

\[X =  \left(X_1\setminus (D^3\times S^1)\right) \cup_{S^2\times
S^1} \left( N\setminus (D^3\times S^1)\right) \]

\noindent where the orbit space of $N$ contains all of the
weighted arcs and isolated fixed points.  The quotient of $X_1$
contains no fixed points but it could still have complicated
topology. $N$ is a simply connected 4--manifold with an $S^1$
action and quotient $S^3$, so it is diffeomorphic to a connect sum
of $S^4$'s, $\CP^2$'s, $\overline{\CP}^2$'s, and $S^2 \times
S^2$'s by Theorem \ref{thm:Fintushel}. We can in fact build $N$ by
starting with a circle action on $S^4$ with two fixed points and
equivariantly connect summing the other factors. This may be yet a
different $S^1$--action, but the resulting manifold is still
diffeomorphic to $N$.  Since $N$ is simply connected, any two
embedded circles are isotopic; hence the fiber sum of $X_1$ with
$N$ along $S^2\times S^1$ using the new $S^1$--action will still
be diffeomorphic to $X$. Once again we can eliminate all connect
sum factors of $N$ except for the $S^4$ we started with.

Thus we have reduced the problem to finding an essential sphere of
nonnegative square in
\[X = \left(X_1\setminus (D^3\times S^1)\right)\cup_{S^2\times
S^1} \left( S^4\setminus (D^3\times S^1)\right) \] where the
quotient of $S^4$ has 2 fixed points. Note that $X$ is a
4--manifold with an $S^1$--action, $b_{+}>0$, and 2 isolated fixed
points $F \subset Y$. Because the sum of the indices is zero, one
of the fixed points comes with a $^+1$ index and the other comes
with a $^-1$ index.  In this situation, $b_+(X)>0$ forces
$b_1(X)>0$ by the formula  $\chi(F) = \chi(X)= 2-2b_1(X) +b_2(X)$
derived from the Smith--Gysin sequence. By Corollary 10.3 of
\cite{circact:Class_of_circ_4man} we have that $b_1(Y)=b_1(X)>0$.
Hence $H_1(Y;\BZ)$ in the long exact sequence of the pair $(Y,F)$
is nonzero:
$$ 0 \ra H_1(Y;\BZ) \ra H_1(Y, F;\BZ) \stackrel{\partial}{\ra}
H_0(F;\BZ) \ra H_0(Y;\BZ) \ra 0.$$ Since $\partial \chi = (1,-1)
\in H_0(F;\BZ)= \BZ \oplus \BZ$ there exists a 1--cycle in
$H_1(X,F;\BZ)$ represented by a closed loop $\gamma$ which is not
a multiple of the Euler class $\chi$ of the action. The preimage
$\pi^{-1}(\gamma)$ is an essential torus of self--intersection 0.
Isotope $\gamma$ to the fixed point set such that $\gamma =
\gamma_1 + \gamma_2$ where $\gamma_i$'s are two arcs running from
one fixed point to the other.  The preimage of both of these arcs
is a sphere of self--intersection 0.  Since the preimage of
$\gamma$ is essential, one of these spheres must be essential as
Theorem \ref{lem:main} demands.

\end{proof}

The following lemma is a slight generalization of Proposition 4 in
\cite{sw:circleactions}.

\begin{lemma} Let $X$ be a smooth, closed oriented $b_+{>}1$ 4--manifold with a
smooth circle action whose orbit space $Y$ has weighted circles
and arcs, isolated fixed points, and one boundary component. Then
there exist an essential sphere of nonnegative square in $X$.
\label{lem:fpbdry}
\end{lemma}

\begin{proof}  As before, we eliminate cases which have spheres of nonnegative square by using
the local description of equivariant plumbing given in Section 4
of \cite{circact:Circ_act_on_four_man}, and by using techniques
used in the proof above.  Thus we can assume that the quotient
space of $X$ contains only simply--weighted circles, $n$ isolated
fixed points $\{x_1, x_2, \ldots, x_n\}$ each with a $^+1$ index,
and one boundary component with index $-n$.  Denote the fixed
point set by $F$.

A linearly independent subset of $H_2(X;\BZ)$ can be constructed
as follows (c.f. Section 8 of
\cite{circact:Circ_act_on_four_man}). For $i=1$ to $n$, let
$\gamma_i$ be an arc that runs from $x_i$ to a point on $\partial
Y$ such that all arcs are mutually disjoint. The preimage
$\pi^{-1}(\gamma_i)$ of each of these arcs is an essential sphere
$S_i$ which represents a 2--cycle in $H_2(X;\BZ)$. These linearly
independent classes have an intersection matrix  with respect to
each other given by
$$S_i \cdot S_j = \left\{ \begin{array}{rll} -1 & \ \ \ & i = j  \\
0 && \mbox{otherwise.} \end{array}\right.$$ Because the
intersection form of $X$ is not negative definite, we must have
$b_2(X)>n$. Let $g$ be the genus of $\partial Y$. Using the fact
$\chi(X)=\chi(F)$, we get
$$b_1(X) = \frac12(b_2(X) - n + 2g).$$
In particular, $b_1(Y) > g$ by Corollary 10.3 of
\cite{circact:Class_of_circ_4man} again.  The long exact sequence
for the pair $(Y,F)$ implies that there is a 1--cycle in
$H_1(Y,F;\BZ)$ which is represented by a  closed loop $\gamma$
which is not a multiple of the Euler class $\chi$. The preimage
$\pi^{-1}(\gamma)$ is an essential torus of self--intersection
zero.  The loop $\gamma$ is homologous to an arc which starts and
ends on $\partial Y$ but is otherwise disjoint from $F$; and the
preimage of the arc is a sphere which is homologous to the torus.
This is an essential sphere of nonnegative square, proving Lemma
\ref{lem:fpbdry}.
\end{proof}






\begin{thebibliography}{99999}

\bibitem{symp:circact:4man_admit_mom_map} K. Ahara, A. Hattori, {\em
Four dimensional symplectic $S^1$--manifolds admitting moment
map}, J. Fac. Sci. Univ. Tokyo, Sect. IA, Math, {\bf 38} (1991),
251 -- 298.

\bibitem{symp:circact:4man_ham_symp} M. Audin, {\em
Hamiltoniens periodiques sur les varietes symplectiques compactes
de dimension 4}, Geometrie symplectique et mechanique, C. Albert,
ed. Lecture Notes in Math. {\bf 1416} (1990).

\bibitem{symp:circact:top_of_torus_act} M. Audin,
`The Topology of Torus Actions on Symplectic Manifolds', Progress
in Mathematics {\bf 93}, Birkh\"{a}user, Basel, 1991.

\bibitem{sw:circleactions} S. Baldridge, {\em Seiberg--Witten invariants of
$4$-manifolds with free circle actions}, Commun. Contemp. Math,
{\bf 3} (2001), 341 -- 353.

\bibitem{sw:sworbcirc} S. Baldridge, {\em Seiberg--Witten invariants, orbifolds, and
circle actions}, E-print GT-0107092.

\bibitem{symp:symp_lef_fib_S1xM} W. Chen and R. Matveyev, {\em Symplectic
Lefschetz fibrations on $S^1\times M^3$}, Geom. Topol. {\bf 4}
(2000), 517--535.

\bibitem{sw:swand4man} S. Donaldson, {\em The Seiberg--Witten equations
and 4-manifold topology}, Bull. A.M.S. {\bf 33} (1996), 45--70.

\bibitem{symp:lefschetz_cpx_strct_seifert}T. Etg\"{u}, {\em  Lefschetz fibrations, complex structures and Seifert fibrations on
$S\sp 1\times M\sp 3$}, Algebr. Geom. Topol. {\bf 1} (2001),
469--489 (electronic).

\bibitem{circact:Circ_act_on_four_man} R. Fintushel, {\em Circle actions
on simply connected $4$--manifolds}, Trans. Amer. Math. Soc. {\bf
230} (1977), 147--171.

\bibitem{circact:Class_of_circ_4man} R. Fintushel,  {\em Classification of
circle actions on 4--manifolds}, Trans. Amer. Math. Soc. {\bf 242}
(1978), 377--390.

\bibitem{sw:class_turkey} R. Fintushel and R. Stern, {\em Immersed
spheres in 4--manifolds and the immersed Thom conjecture}, Turkish
J. Math. {\bf 19} (1995), 145--157.

\bibitem{sw:4man_kirby_calc} R. Gompf and A. Stipsicz,
`4--Manifolds and Kirby Calculus', Graduate Studies in Mathematics
{\bf 20}, American Mathematical Society, Providence, Rhode Island,
1999.


\bibitem{sw:symp:embedded_spheres} T.J. Li, {\em Smoothly embedded spheres in symplectic
4--manifolds}, Proc. Amer. Math. Soc.  {\bf 127}  (1999), no. 2,
609--613.

\bibitem{symp:immersed_spheres} D. McDuff, {\em Immersed spheres in symplectic 4--manifolds},
Ann. Inst. Fourier, Grenoble {\bf 42} (1990), 369-392.

\bibitem{symp:momentmap} D. McDuff, {\em The moment map for circle actions
on symplectic manifolds}, J. Geom. Phys.  {\bf 5} (1988),
149--160.

\bibitem{sw:sw_inv_seifert_space} T. Mrowka, P. Ozsv\'{a}th, and B. Yu, {\em
Seiberg--Witten monopoles on Seifert fibered spaces}, Comm. Anal.
Geom. {bf 5} (1997), 685 -- 791.

\bibitem{sw:abiabatic_limits} L. Nicolaescu, {\em Adiabatic limits of the Seiberg--Witten
equations on Seifert manifolds}, Comm. in Anal. and Geom. {\bf 6}
(1998), 61-123.

\bibitem{sw:notes_on_sw_theory} L. Nicolaescu, 'Notes on Seiberg--Witten Theory',
Graduate Studies in Mathematics {\bf 28}, American Mathematical
Society, Providence, Rhode Island, 2000.

\bibitem{circact:nonlincircact} P. Pao, {\em Nonlinear circle actions on the
$4$--sphere and twisting spun knots.}  Topology {\bf 17} (1978),
no. 3, 291--296.

\bibitem{symp:circact:semifree_symp_circ_act} S. Toleman and J. Weitsman,
{\em On semifree symplectic circle actions with isolated fixed
points}, Topology {\bf 39} (2000), 299--309.

\bibitem{sw:sw_inv_and_symp_form} C. Taubes, {\em The Seiberg--Witten
invariants and symplectic forms}, Math. Res. Lett.   {\bf 1}
(1994), no. 6, 809--822.

\bibitem{sw:geometry_of_sw} C. Taubes, {\em The geometry of the Seiberg--Witten invariants},
Surveys in differential geometry, Vol. III (Cambridge, MA, 1996)
299--339, Int. Press, Boston, MA, 1998.


\bibitem{sw:monopolesandfourman} E. Witten, {\em Monopoles and
four--manifolds}, Math. Res. Lett. {\bf 1} (1994), no. 6,
769--796.

\end{thebibliography}
\end{document}